\documentclass[letterpaper, 10pt, journal, twocolumn, final]{IEEEtran}

\usepackage{amsthm}
\usepackage{soul}
\usepackage[normalem]{ulem}
\usepackage[hidelinks]{hyperref}
\usepackage{bookmark}
\usepackage[latin1]{inputenc}
\usepackage{microtype}
\usepackage[noadjust]{cite}

\usepackage[dvipsnames,table]{xcolor}
\usepackage{amsfonts,amssymb,amsmath,color,amsthm}
\usepackage{placeins}
\usepackage{graphicx}
\usepackage{algorithm}
\usepackage[]{algpseudocode}

\usepackage{tikz}
\usetikzlibrary{shapes,arrows}

\usepackage{array}
\usepackage{booktabs}

\usepackage{dsfont} %

\newcommand{\map}[3]{#1: #2 \rightarrow #3}

\newcommand{\EE}{\mathcal{E}} 
\newcommand{\GG}{\mathcal{G}}

\newcommand{\subj}{\textnormal{subj. to}}

\newcommand{\nbrs}{\mathcal{N}}

\newcommand{\real}{{\mathbb{R}}}
\renewcommand{\natural}{{\mathbb{N}}}
\newcommand{\integer}{{\mathbb{Z}}}
\newcommand{\conv}[1]{\textnormal{conv}(#1)}

\newcommand{\until}[1]{\{1,\ldots,#1\}} 
\newcommand{\fromto}[2]{\{#1,\ldots,#2\}}

\newcommand\oprocendsymbol{\hbox{$\square$}}
\newcommand\oprocend{\relax\ifmmode\else\unskip\hfill\fi\oprocendsymbol}

\newtheorem{theorem}{Theorem}[section]

 \newtheorem{lemma}[theorem]{Lemma}
\newtheorem{remark}[theorem]{Remark}

\newtheorem{assumption}[theorem]{Assumption}

\newcommand{\0}{0}
\newcommand{\1}{\mathds{1}}

\newcommand{\bx}{x}
\newcommand{\bz}{z}

\newcommand{\by}{y}

\newcommand{\bmu}{\mu}

\newcommand{\kron}{\otimes}
\newcommand{\smallsum}{\textstyle\sum\limits}

\newcommand{\algx}{\bx^\infty}

\renewcommand{\Xi}{X_i^{\textsc{milp}}}

\newcommand{\dz}{p}
\newcommand{\dr}{q}

\newcommand{\INT}{I_{\mathbb{Z}}}

\newcommand{\LP}{\textsc{lp}}

\newcommand{\bL}{\ell}

\newcommand{\agents}{\mathbb{I}}

\newcommand{\GEN}{\agents_\textsc{gen}}
\newcommand{\CLOAD}{\agents_\textsc{cl}}
\newcommand{\LOAD}{\agents_\textsc{lo}}
\newcommand{\STOR}{\agents_\textsc{stor}}
\newcommand{\UGRID}{\agents_\textsc{grid}}
\newcommand{\REN}{\agents_\textsc{ren}}
\newcommand{\pl}{\textsc{pl}}
\newcommand{\bE}{E}

\newcommand{\xmin}{x^\textsc{min}}
\newcommand{\xmax}{x^\textsc{max}}
\newcommand{\OM}{\zeta}

\newcommand{\expv}{{\mathbb{E}}}

\renewcommand{\Xi}{X_i^{\textsc{milp}}}
\newcommand{\convXi}{\conv{\Xi}}

\newcommand{\interv}[1]{1,\ldots,#1}

\newcommand{\xl}{\bx^{\textsc{L}}}
\newcommand{\etal}{\eta^{\textsc{L}}}

\renewcommand{\algx}{\bx^\infty}
\newcommand{\algeta}{\eta^\infty}

\makeatletter
\newcommand{\StatexIndent}[1][3]{%
  \setlength\@tempdima{\algorithmicindent}%
  \Statex\hskip\dimexpr#1\@tempdima\relax}
\makeatother

\renewcommand{\lim}{\operatornamewithlimits{lim\vphantom{p}}}

\graphicspath{{figs/},{figs/simulations/}}

\def \sgalgname/{Distributed Stochastic Mixed-integer Microgrid Control}

\begin{document}

\title{A Distributed Mixed-Integer Framework to \\ Stochastic Optimal Microgrid Control}

\author{Andrea Camisa, %
  Giuseppe Notarstefano %
  \thanks{A. Camisa and G. Notarstefano are with the Department of Electrical, 
  Electronic and Information Engineering, University of Bologna, Bologna, Italy. 
  \texttt{\{a.camisa, giuseppe.notarstefano\}@unibo.it}.
  This result is part of a project that has received funding from the European 
  Research Council (ERC) under the European Union's Horizon 2020 research 
  and innovation programme (grant agreement No 638992 - OPT4SMART).
  }
}

\maketitle

\begin{abstract}
  This paper deals with distributed control of microgrids composed
  of storages, generators, renewable energy sources, critical and
  controllable loads.
  We consider a stochastic formulation of the optimal control problem
  associated to the microgrid that appropriately takes into account
  the unpredictable nature of the power generated by renewables.
  The resulting problem is a Mixed-Integer Linear Program and is
  NP-hard and nonconvex.
  Moreover, the peculiarity of the considered framework is that no
  central unit can be used to perform the optimization, but rather
  the units must cooperate with each other by means of neighboring
  communication.
  To solve the problem,
  we resort to a distributed methodology based on a primal decomposition approach.
  The resulting algorithm is able to compute high-quality feasible solutions
  to a two-stage stochastic optimization problem, for which we also provide a theoretical
  upper bound on the constraint violation.
  Finally, a Monte Carlo numerical computation on a scenario with
  a large number of devices shows the efficacy of the
  proposed distributed control approach.
  The numerical experiments are performed on realistic scenarios
  obtained from Generative Adversarial Networks trained
  an open-source historical dataset of the EU.
\end{abstract}

\section{Introduction}
\label{sec:intro}
In the last decade, the use of renewable energy sources is soaring and
is creating new challenges in the field of microgrid control. These important
structural changes of the power grid call for novel approaches that
must appropriately take into account the stochastic nature of the energy
produced by renewables.
To this end, optimization-based control techniques are increasingly used.
However they typically employ centralized approaches that require the
collection of the problem data at each node, which may lead to
a single point of failure. Distributed optimization approaches are a
promising alternative that allows for the solution of optimization
problems with spatially distributed data while preserving the locality
of the data and even resilience of the network in case of
failures~\cite{molzahn2017survey,nedic2018distributed,notarstefano2019distributed}.
We first review optimal control techniques, then we
recall approaches based on mixed-integer programming and finally
move to distributed approaches.
Optimal control techniques allow for shaping input trajectories
that take into account energy consumption/production costs and
user comfort. In recent times, they are increasingly achieved with
moving horizon techniques as Model Predictive Control (MPC)
as it flexibly allows one to tackle several challenges,
see e.g.~\cite{ouammi2015coordinated,cominesi2017two,le2017plug}.
Stochastic optimization-based approaches are also being developed.
In~\cite{zakariazadeh2014smart},
a stochastic optimization method for energy and reserve scheduling
with renewable energy sources and demand-side participation is considered.
The work~\cite{nguyen2015stochastic} studies a stochastic unit commitment
and economic dispatch problem with renewables and incorporating the
battery operating cost.
Another prominent approach is Mixed-Integer Linear
Programming (MILP), which is gathering significant attention
due to its ability to model logical statements that often occur
within microgrids.
In~\cite{kriett2012optimal}, a MILP optimal control approach of residential microgrid
is proposed.
In~\cite{marzband2014experimental} a mixed-integer nonlinear programming
formulation is considered with experimental validation for islanded-mode microgrids.
In~\cite{shirazi2017cost}, a MILP is formulated to achieve optimal load shifting in
microgrids.
The MPC and the MILP approaches have been combined in~\cite{parisio2014model},
which proposes a receding horizon implementation of the MILP approach on an
experimental testbed. A stochastic version of this work is considered in~\cite{parisio2016stochastic},
which further takes into account renewable energy sources and
aims at an environmental/economical operation of microgrids.
While these works take into account more and more aspects
of microgrids, they are all based on centralized optimization techniques
that require one of the nodes to be chosen as master, thus introducing
scalability and privacy issues.
As energy networks are intrinsically distributed, there is often the need to
devise distributed approaches %
that exploit the graph structure.
The work~\cite{molzahn2017survey} reviews distributed
methods for optimal power flow problems, while~\cite{dorfler2019distributed}
surveys distributed control approaches for autonomous
power grids.
In~\cite{bolognani2013distributed}, a distributed approach to optimal
reactive power compensation is proposed. 
In~\cite{causevic2018energy,belluschi2020distributed}, the authors
propose distributed algorithms for optimal energy building management,
while~\cite{cavraro2020distributed} investigates a distributed feedback
control law to minimize power generation cost in prosumer-based
networks.
However, none of the mentioned works formulates a comprehensive
stochastic scheduling problem involving the demand-side in a distributed way.
Novel distributed methods relying on MILPs can take advantage of the
latest progress of distributed optimization methods.
MILPs are nonconvex and NP-hard, therefore large-scale instances
can be solved within acceptable time windows only suboptimally.
On this regard, the recent works~\cite{falsone2018distributed,camisa2021distributed}
propose distributed algorithms to compute feasible solutions of MILPs over networks.

The contributions of this paper are as follows. We consider a distributed stochastic microgrid control problem
consisting of several interconnected power units, namely generators, renewable
energy sources, storages and loads. We begin by recalling the microgrid model.
We then show that the optimal control problem
can be recast as a distributed MILP.
We then apply a two-stage stochastic programming approach
to the distributed MILP and show that also this problem can be
cast as a distributed MILP. This new problem is then tackled using
an approach inspired to recent approaches proposed in the literature, %
which are suitably modified to deal with the stochastic scenario.
The proposed algorithm provides a feasible solution to the two-stage
stochastic problem at each iteration, while preserving sensible data at each
node. As the algorithm progresses, the cost of the provided
solution improves and the expected violation of the power balance
constraint decreases.
For the asymptotic solution provided by the algorithm, we formally
prove an upper bound on the violation of the power balance
constraint.
  We then apply the developed approach to a simulation scenario with
  a large number of devices.
  We perform realistic simulations by using open-source historical data,
  taken from the EU platform Open Power System Data~\cite{opsd2020timeseries},
  on energy generation/consumption in South Italy. We train a
  Generative Adversarial Neural Network (GAN) based on these data and use
  it to generate sample energy generation/consumption profiles.
  The generated data is used to perform a Monte Carlo numerical
  experiment on the Italian HPC CINECA infrastracture
  to show the efficacy of the distributed algorithm.

The paper is organized as follows. In Section~\ref{sec:microgrid_deterministic}, we
describe the mixed-integer microgrid model and the stochastic
optimal control problem.
In Section~\ref{sec:distributed_stochastic}, we reformulate the problem as a
distributed MILP and apply the two-stage stochastic programming approach.
In Section~\ref{sec:algorithm}, we describe the proposed distributed algorithm
and provide theoretical results on the worst-case constraint violation,
while in Section~\ref{sec:simulations} we discuss Monte Carlo numerical simulations
on a practical scenario with a large number of devices and realistic synthesized data.

\section{Stochastic Mixed-Integer Microgrid Control with Renewables}
\label{sec:microgrid_deterministic}
\begin{table}[t]\centering
	\caption{List of the main symbols and their definitions}
	\label{tb:symbols}
	{\rowcolors{2}{gray!20}{white}
	\renewcommand{\arraystretch}{1.15}
	\begin{tabular}{ll}
	  \rowcolor{gray!30}
	  \hline
		\multicolumn{2}{c}{Basic definitions}
		\\
		\hline
		$N\in\mathbb{N}$ 
		& Number of units in the system
		\\
		$\agents = \until{N}$
		& Set of units
		\\
		$\varepsilon>0$
		& Very small number (e.g. machine precision)
		\\
		\rowcolor{gray!30}
		\hline
		\multicolumn{2}{c}{Storages (indexed by $i \in \STOR$)}
		\\
		\hline
		$x_i(k)$
		& State of charge at time $k$
		\\
		$u_i(k)$
		& Exchanged power ($\ge 0$ if charging) at time $k$
		\\
		$\delta_i(k)$
		& Charging (1) / discharging (0) state
		\\
		$z_i(k)$
		& Auxiliary optimization variable
	  \\
		$\eta_i^c$, $\eta_i^d$
		& Charging and discharging efficiencies
	  \\
		$\xmin_i$, $\xmax_i$
		& Minimum and maximum storage level
		\\
		$x_i^\pl$
		& Physiological loss of energy
		\\
		$C_i$
		& Maximum output power
		\\
		$\OM_i$
		& Operation and maintenance cost coefficient
		\\
		\rowcolor{gray!30}
		\hline
		\multicolumn{2}{c}{Generators (indexed by $i \in \GEN$)}
		\\
		\hline
		$u_i(k)$
		& Generated power ($\ge 0$) at time $k$
		\\
		$\delta_i(k)$
		& On (1) / off (0) state (``on'' iff $u_i(k) > 0$)
		\\
		$\nu_i(k)$
		& Epigraph variable for quadratic generation cost
		\\
		$\theta_i^\textsc{u}(k)$, $\theta_i^\textsc{d}(k)$
		& Epigraph variables for startup/shutdown costs
		\\
		$T_i^\textsc{up}$, $T_i^\textsc{down}$
		& Minimum up/down time
		\\
		$u_i^\textsc{min}$, $u_i^\textsc{max}$
		& Min. and max. power that can be generated
		\\
		$r_i^\textsc{max}$
		& Maximum ramp-up/ramp-down
		\\
		$\kappa_i^\textsc{u}(k)$, $\kappa_i^\textsc{d}(k)$
		& Startup and shutdown costs
		\\
		$\OM_i$
		& Operation and maintenance cost coefficient
		\\
		\rowcolor{gray!30}
		\hline
	  \multicolumn{2}{c}{Renewable energy sources (indexed by $i \in \REN$)}
		\\
		\hline
		$P_i(k)$
		& Generated power at time $k$
		\\
		\rowcolor{gray!30}
		\hline
	  \multicolumn{2}{c}{Controllable loads (indexed by $i \in \CLOAD$)}
		\\
		\hline
		$\beta_i(k)$
		& Curtailment factor ($\in [\beta_i^\textsc{min}, \beta_i^\textsc{max}]$) at time $k$
		\\
		$D_i(k)$
		& Consumption forecast at time $k$
		\\
		$\beta_i^\textsc{min}$, $\beta_i^\textsc{max}$
		& Minimum and maximum allowed curtailment
    \\
    \rowcolor{gray!30}
		\hline
	  \multicolumn{2}{c}{Connection to the main grid (indexed by $i \in \UGRID$)}
		\\
		\hline
		$u_i(k)$
		& Imported power from the grid at time $k$
		\\
		$\delta_i(k)$
		& Importing (1) or exporting (0) mode at time $k$
		\\
		$\phi_i(k)$
		& Total expenditure for imported power at time $k$
		\\
		$\phi_i^\textsc{p}(k)$, $\phi_i^\textsc{s}(k)$
		& Price for power purchase and sell at time $k$
		\\
		$P_i^\textsc{max}$
		& Maximum exchangeable power
		\\
		\rowcolor{gray!30}
		\hline
	  \multicolumn{2}{c}{Two-stage stochastic problem}
		\\
		\hline
		$q_+, q_-$
		& Costs associated to positive and negative recourse
	\end{tabular}}
\end{table}

Let us begin by introducing the mixed-integer microgrid model.
For ease of exposition, we consider a fairly general model
inspired to the one in~\cite{parisio2016stochastic}
without taking into account some specific aspects
(see also Remark~\ref{rem:chp}).
This allows us to better highlight the main features of
the proposed approach while keeping the discussion not
too technical.
A microgrid consists of $N$ units, partitioned as follows.
Storages are collected in $\STOR$, generators in $\GEN$,
renewable energy sources in $\REN$,
critical loads in $\LOAD$, controllable loads in $\CLOAD$
and one connection with the utility grid in $\UGRID$.
The whole set of units is
\begin{align*}
  \agents = \until{N} = \STOR \cup \GEN \cup \REN \cup \LOAD \cup \CLOAD \cup \UGRID.
\end{align*}
Throughout the document, we interchangeably refer to the
units also as \emph{agents}.
In the next subsections we describe each type of units separately, while in
Section~\ref{sec:opt_control_pb} we will introduce the optimal control problem.
In the following, we denote the optimal control prediction horizon as $T \in \natural$.

\subsection{Storages}
For storage units $i \in \STOR$, let $x_i(k) \in \real$ be the stored energy level at time
$k$ and let $u_i(k) \in \real$ denote the power exchanged with the storage unit at time $k$
(positive for charging, negative for discharging).
The dynamics at each time $k$ amounts to $x_i(k+1) = x_i(k) + \eta_i u_i(k) - x_i^\pl$,
where $\eta_i$ denotes the (dis)charging efficiency and $x_i^\pl$ is a
physiological loss of energy. It is assumed that $\eta_i = \eta_i^c$ if
$u_i(k) \ge 0$ (charging mode), whereas $\eta_i = 1/\eta_i^d$ if $u_i(k) < 0$
(discharging mode), with $0 < \eta_i^c, \eta_i^d < 1$.
Thus, the dynamics is piece-wise linear. To deal with this,
we utilize mixed-integer
inequalities~\cite{bemporad1999control}. Let us
introduce additional variables $\delta_i(k) \in \{0,1\}$ and
$z_i(k) \triangleq \delta_i(k) u_i(k) \in \real$ for all $k$.
Each $\delta_i(k)$ is one if and only if $u_i(k) \ge 0$
(i.e. the storage unit at time $k$ is in the charging state).
After following the manipulations proposed in~\cite{parisio2014model},
we obtain the following model for the $i$-th storage unit,
\begin{subequations}
\begin{align}
  &x_i(k+1) = x_i(k) + (\eta_i^c - \tfrac{1}{\eta_i^d}) z_i(k) + \tfrac{1}{\eta_i^d} u_i(k) - x_i^\pl,
	\label{eq:storage_1}
	\\
  &\bE_i^1 \delta_i(k) + \bE_i^2 z_i(k) \le \bE_i^3 u_i(k) + \bE_i^4,
	\label{eq:storage_2}
	\\
  &x_i^\textsc{min} \le x_i(k) \le x_i^\textsc{max},
	\label{eq:storage_3}
\end{align}
for all time instants $k$, and %
\begin{align}
  x_i(0) = x_{i,0},
  \label{eq:storage_4}
\end{align}
\label{eq:storage}%
\end{subequations}
where~\eqref{eq:storage_1} is the dynamics,~\eqref{eq:storage_2} are
mixed-integer inequalities expressing the logical constraints,~\eqref{eq:storage_3}
are box constraints on the
state of charge (with $0 < x_i^\textsc{min} < x_i^\textsc{max}$),
and~\eqref{eq:storage_4} imposes the initial condition ($x_{i,0} \in \real$
is the initial state of charge of storage $i$). The matrices in~\eqref{eq:storage_2} are
\begin{align*}
	\bE_i^1 \!=\!\! \begin{bmatrix}C_i \\ -(C_i \!+\! \varepsilon) \\ C_i \\ C_i \\ -C_i \\ -C_i\end{bmatrix}\!, \:\:
	\bE_i^2 \!=\!\! \begin{bmatrix} 0 \\ 0 \\ 1 \\ -1 \\ 1 \\ -1\end{bmatrix}\!, \:\:
	\bE_i^3 \!=\!\! \begin{bmatrix} 1 \\ -1 \\ 1 \\ -1 \\ 0 \\ 0\end{bmatrix}\!, \:\:
	\bE_i^4 \!=\!\! \begin{bmatrix} C_i \\ -\varepsilon \\ C_i \\ C_i \\ 0 \\ 0\end{bmatrix}\!,
\end{align*}
where $C_i > 0$ is the limit output power and
$\varepsilon > 0$ is a very small number (typically machine precision).
To each storage $i$ it is associated an operation and maintenance cost,
which is equal to
\begin{align}
  J_i = \sum_{k=0}^{K-1} \OM_i |u_i(k)|
  = \sum_{k=0}^{K-1} \OM_i (2 z_i(k) - u_i(k)),
\end{align}
where $\OM_i > 0$ is the operation and maintenance cost per exchanged unit of power
and $2 z_i(k) - u_i(k) = |u_i(k)|$ is the absolute value of the power exchanged
with the storage.

\subsection{Generators}
\begin{subequations}
\label{eq:generator}
For generators $i \in \GEN$, let $u_i(k) \in \real, u_i(k) \ge 0$ denote the
generated power at time $k$. Since generators can be either on or off, as
done for the storages we let $\delta_i(k) \in \{0,1\}$ be an
auxiliary variable that is equal to $1$ if and only if $u_i(k) > 0$.
As in the case of storages, we must consider constraints on the operating
conditions of generators. Namely, if a generator is turned on/off, there is
a minimum amount of time for which the unit must be kept on/off.
This logical constraint is modeled by the inequalities
\begin{align}
  &\delta_i(k) - \delta_i(k-1) \le \delta_i(\tau),
  \nonumber
  \\
  & \hspace{2cm}
  \tau = k+1, \ldots, \min(k+T_i^\textsc{up} - 1, T),
  \label{eq:generator_1}
\end{align}
\begin{align}
  &\delta_i(k-1) - \delta_i(k) \le \delta_i(\tau),
  \nonumber
  \\
  & \hspace{2cm}
  \tau = k+1, \ldots, \min(k+T_i^\textsc{down} - 1, T),
  \label{eq:generator_2}
\end{align}
for all time instants $k$, %
where $T_i^\textsc{up}$ and
$T_i^\textsc{down}$ are the minimum up and down time of
generator $i$.
The power flow limit and the ramp-up/ramp-down limits are modeled
respectively by
\begin{align}
  u_i^\textsc{min} \delta_i(k) &\le u_i(k) \le u_i^\textsc{max} \delta_i(k),
  \label{eq:generator_3}
  \\
  -r_i^\textsc{max} \delta_i(k) &\le u_i(k) - u_i(k-1) \le r_i^\textsc{max} \delta_i(k),
  \label{eq:generator_4}
\end{align}
for all times $k$, %
where $u_i^\textsc{max} \ge u_i^\textsc{min} \ge 0$ denote
the maximum and minimum power that can be generated by generator $i$
and $r_i^\textsc{max} \ge 0$ denotes the maximum ramp-up/ramp-down.

The cost associated to generator units is composed of three parts,
which are \emph{(i)} a (quadratic) generation cost,
\emph{(ii)} a start-up/shut-down cost and
\emph{(iii)} an operation and maintenance cost.
To model the generation cost, %
we
consider a piece-wise linearized version
$\max\limits_{\ell} \big( S_i^\ell u_i(k) + s_i^\ell \big)$
for all $k$ with appropriately defined $S_i^\ell, s_i^\ell \in \real$.
The startup $\theta_i^\textsc{u}$ and shutdown cost $\theta_i^\textsc{d}$ at each time
$k \in \fromto{0}{K-1}$ are equal to
\begin{align*}
  \theta_i^\textsc{u}(k) &= \max\Big\{ 0, \:\: \kappa_i^\textsc{u}(k)[\delta_i(k) - \delta_i(k-1)] \Big\},
  \\
  \theta_i^\textsc{d}(k) &= \max\Big\{ 0, \:\: \kappa_i^\textsc{d}(k)[\delta_i(k-1) - \delta_i(k)] \Big\},
\end{align*}
where $\kappa_i^\textsc{u}(k), \kappa_i^\textsc{d}(k) > 0$
are the start-up and shut-down costs at time $k$.
The operation and maintenance cost is equal to $\zeta_i \delta_i(k)$,
where $\zeta_i > 0$ is a cost coefficient (we assume that there is no cost
when the generator is turned off).
Thus, the expression for the cost of each generator $i$ is
\begin{align*}
  J_i =
  \!\sum_{k=0}^{K-1} \!\Big[
    \max\limits_{\ell} \big( S_i^\ell u_i(k) + s_i^\ell \big) + \zeta_i \delta_i(k)
    + \theta_i^\textsc{u}(k) + \theta_i^\textsc{d}(k)
  \Big].
\end{align*}
Note that the cost function has internal maximizations and,
as such, is nonlinear.
However, since the cost is to be minimized, it can be recast as a linear function
by introducing so-called epigraph variables (see e.g.~\cite{boyd2004convex})
as follows.
As regards the generation cost, we replace it
with epigraph variables $\nu_i(k) \in \real$ and impose the constraints
\begin{align}
  \nu_i(k) \ge S_i^\ell u_i(k) + s_i^\ell, \hspace{1cm} \forall \: \ell,
\label{eq:generator_5}
\end{align}
for all times $k$. %
Similarly, we can treat $\theta_i^\textsc{u}, \theta_i^\textsc{d} \in \real$
as epigraph variables and write the constraints
\begin{align}
  \theta_i^\textsc{u}(k) &\ge \kappa_i^\textsc{u}(k)[\delta_i(k) - \delta_i(k-1)],
  \label{eq:generator_6}
  \\
  \theta_i^\textsc{d}(k) &\ge \kappa_i^\textsc{d}(k)[\delta_i(k-1) - \delta_i(k)],
  \label{eq:generator_7}
  \\
  \theta_i^\textsc{u}(k) &\ge 0,
  \label{eq:generator_8}
  \\
  \theta_i^\textsc{d}(k) &\ge 0,
  \label{eq:generator_9}
\end{align}\end{subequations}
for all $k$. We therefore obtain the following expression
for the cost function of generator $i$,
\begin{align}
  J_i = \sum_{k=0}^{K-1} \big[ \nu_i(k) + \theta_i^\textsc{u}(k) + \theta_i^\textsc{d}(k) + \zeta_i \delta_i(k) \big].
\end{align}

\subsection{Renewable Energy Sources}

We consider two types of renewables, namely wind generators
and solar generators.
Rather than using a physical or dynamical model for these generators,
we use a predictor to generate realistic power production scenarios.
Indeed, thanks to the huge amount of historical datasets freely available
on the internet, neural network-based predictors have excellent
accuracy. More details are in Section~\ref{sec:gan}.
We will employ this technique also to generate power demand predictions.%
These units only contribute to the power balance
constraint~\eqref{eq:power_balance} through their
generated power at each time slot $k$, denoted as $P_i(k) \ge 0$,
and do not have associated cost or constraints.
Note that $P_i(k)$ are unknown beforehand and must be modeled as
stochastic variables having a certain probability distribution. We
discuss this aspect more in details in Section~\ref{sec:distributed_stochastic}.

\subsection{Loads}
We consider two types of loads, namely critical loads and controllable loads.
For critical loads $i \in \LOAD$, we will denote by $D_i(k)$ the consumption
forecast at time $k$ and we assume it is given. There are no
optimization variables (and thus cost functions) associated with this kind
of units, however their consumption must be considered in the power balance
(cf. Section~\ref{sec:opt_control_pb}). %

For controllable loads $i \in \CLOAD$, let $D_i(k)$ be the consumption
forecast at time $k$, which is assumed to be given.
In case the microgrid has energy shortages, the consumption
of controllable loads can be curtailed to meet power balance
constraints. This is quantified with a curtailment factor
$\beta_i(k) \in [\beta_i^\textsc{min},\beta_i^\textsc{max}]$,
where $0 \le \beta_i^\textsc{min} \le \beta_i^\textsc{max} \le 1$ are
the bounds on the allowed curtailment.
The actual power consumption at time $k$
is thus $(1 - \beta_i(k)) D_i(k)$, i.e. if $\beta_i(k) = 0$ there is
no curtailment. The curtailment factor is an optimization variable
and can be freely chosen, thus in principle it can be $\beta_i(k) >0$
for some $k$ (even if there are no energy shortages) if this results in
a cost improvement. The following constraint must be imposed, %
\begin{align}
  \beta_i^\textsc{min} \le \beta_i(k) \le \beta_i^\textsc{max},
\label{eq:load}
\end{align}
for all times $k$.
We assume the microgrid incurs in a cost that is proportional to the
total curtailed power, thus the cost function associated to controllable load $i$ is
\begin{align}
  J_i = \sum_{k=0}^{K-1} \varphi_i D_i(k) \beta_i(k),
\end{align}
where $\varphi_i > 0$ is a penalty weight.

\subsection{Connection to the Utility Grid}
For the connection with the utility grid $i \in \UGRID$, let $u_i(k) \in \real$
denote the imported (exported) power level from (to) the
utility grid. We use the convention that imported power at time $k$
is non-negative $u_i(k) \ge 0$.
As before, since the power purchase price is different from the power sell price,
we consider auxiliary optimization variables $\delta_i(k) \in \{0,1\}$ and $\phi_i(k) \in \real$ for all $k$.
The variable $\delta_i(k)$ models the logical statement $\delta_i(k) = 1$ if and only
if $u_i(k) \ge 0$ (i.e. power is imported from the utility grid).
The variable $\phi_i(k)$ represents the total expenditure (retribution)
for imported (exported) energy.
Denoting by $\phi_i^\textsc{p}(k), \phi_i^\textsc{s}(k) \ge 0$ the price for power
purchase and sell, it holds $\phi_i(k) = \phi_i^\textsc{p}(k) u_i(k)$ if
$\delta_i(k) = 1$ and $\phi_i(k) = \phi_i^\textsc{s}(k) u_i(k)$ if $\delta_i(k) = 0$.
By denoting by $P_i^\textsc{max} \ge 0$ the maximum exchangeable power,
the corresponding mixed-integer inequalities are (cf.~\cite{parisio2014model}),
\begin{align}
  \bE_i^1 \delta_i(k) + \bE_i^2 \phi_i(k) \le \bE_i^3(k) u_i(k) + \bE_i^4,
\label{eq:grid}
\end{align}
for all $k$,
where the matrices are defined as
\begin{align*}%
	\bE_i^1 \!=\!\!\! \begin{bmatrix} P_i^\textsc{max} \\ \!-\!P_i^\textsc{max} \!\!-\! \varepsilon\!\! \\ M_i \\ M_i \\ -M_i \\ -M_i \end{bmatrix}\!\!, 
	\bE_i^2 \!=\!\!\! \begin{bmatrix} 0 \\ 0 \\ 1 \\ -1 \\ 1 \\ -1 \end{bmatrix}\!\!, 
	\bE_i^3(k) \!=\!\!\! \begin{bmatrix} 1 \\ -1 \\ \phi^\textsc{p}(k)\! \\ \!-\phi^\textsc{p}(k)\!\! \\ \phi^\textsc{s}(k)\! \\ \!-\phi^\textsc{s}(k)\!\! \end{bmatrix}\!\!, 
	\bE_i^4 \!=\!\!\! \begin{bmatrix} \!P_i^\textsc{max}\!\! \\ -\varepsilon \\ M_i \\ M_i \\ 0 \\ 0\end{bmatrix}\!\!,
\end{align*}
with $M_i = P_i^\textsc{max} \cdot \max\limits_k(\phi^\textsc{p}(k), \phi^\textsc{s}(k))$.
Clearly, the cost associated with this unit is
\begin{align}
  J_i = \sum_{k=0}^{K-1} \phi_i(k).
\end{align}

\subsection{Power Balance Constraint and Optimal Control Problem}
\label{sec:opt_control_pb}
Electrical balance must be met at each time $k$, i.e.,
\begin{align}
  &u_{\UGRID}(k)
  = 
  \sum_{i \in \STOR}\! u_i(k) -\! \sum_{i \in \GEN} u_i(k) +\! \sum_{i \in \CLOAD} (1 - \beta_i(k)) D_i(k)
  \nonumber
  \\
  &\hspace{1.6cm}
  + \sum_{i \in \LOAD} D_i(k) - \sum_{i \in \REN}\! P_i(k),
\label{eq:power_balance}
\end{align}
Recall that the length of the prediction horizon is $K \in \natural$.
The optimal control problem, which is a MILP, can be posed as
\begin{align}
  \min_{u} \: & \sum_{k=0}^{K-1} \!\bigg[ \phi_\textsc{grid}(k) \!
    +\!\!\! \sum_{i \in \GEN}\! (\OM_i \delta_i(k) \!+\! \nu_i(k) \!+\! \theta_i^\textsc{u}(k) \!+\! \theta_i^\textsc{d}(k)) 
  \nonumber
  \\
  & \hspace{0.3cm}
    +\sum_{i \in \CLOAD} \varphi_i D_i(k) \beta_i(k) +\!\!\sum_{i \in \STOR} \!\!\big(\OM_i (2 z_i(k) - u_i(k)) \big)\!\bigg]
  \nonumber
  \\
  \subj \: & \: \text{storage constraints } \eqref{eq:storage}
  \label{eq:microgrid_mpc_problem}
  \\
  & \: \text{generator constraints } \eqref{eq:generator}
  \nonumber
  \\
  & \: \text{constraints } \eqref{eq:load}, \eqref{eq:grid}, \eqref{eq:power_balance}.
  \nonumber
\end{align}
Note that problem~\eqref{eq:microgrid_mpc_problem} is a stochastic
optimization problem. Indeed, the equality constraint~\eqref{eq:power_balance}
is stochastic since it depends on $P_i(k)$. Next we show how to handle this
level of complexity.

\begin{remark}
\label{rem:chp}
  Note that the microgrid model can also be extended to additionally
  consider thermal loads and Combined Heat and Power (CHP) units,
  which would additionally require
  a thermal balance constraint. The architecture proposed in the
  following can be easily adapted to deal with this scenario by making
  only minor changes. However, in order not to complicate too much
  the notation, we prefer not to introduce this further level of complexity,
  which nevertheless can be handled by the proposed framework.
  \oprocend
\end{remark}

\section{Distributed Constraint-coupled Stochastic Optimization}
\label{sec:distributed_stochastic}
To handle the stochastic quantities $P_i(k)$ we follow the ideas
of~\cite{parisio2016stochastic} and utilize a two-stage stochastic
optimization approach.
As we are interested in a distributed algorithm, instead of
applying the two-stage stochastic approach directly to problem~\eqref{eq:microgrid_mpc_problem},
we rather apply it to a distributed reformulation of problem~\eqref{eq:microgrid_mpc_problem}.
In this section, we first introduce the distributed reformulation of the
problem and then formalize the two-stage stochastic optimization approach.

\subsection{Constraint-coupled Reformulation}
\label{sec:constraint_coupled_formulation}

The optimal control problem~\eqref{eq:microgrid_mpc_problem} can be
reformulated in such a way that the distributed structure of the problem
becomes more evident. Formally, problem~\eqref{eq:microgrid_mpc_problem}
is equivalent to the stochastic \emph{constraint-coupled} MILP,
\begin{align}
\begin{split}
  \min_{\bx_1,\ldots,\bx_N} \: & \: \sum_{i =1}^N c_i^\top \bx_i
  \\
  \subj \: 
  & \: \sum_{i=1}^N A_i \bx_i = b
  \\
  & \: \bx_i \in \Xi, \hspace{1cm} i = \interv{N},
\end{split}
\label{eq:MILP}
\end{align}
where, for all $i \in \until{N}$, the decision vector $\bx_i$ has $n_i = \dz_i + \dr_i$
components (thus $c_i \in \real^{n_i}$) with $\dz_i, \dr_i \in \natural$
and the local constraint set is of the form
\begin{align*}
  \Xi = P_i \cap (\integer^{\dz_i} \times \real^{\dr_i}),
\end{align*}
for some nonempty compact polyhedron $P_i \subset \real^{\dz_i + \dr_i}$.
Moreover, the matrices $A_i \in \real^{K \times n_i}$
and the vector $b \in \real^K$ model coupling constraints among the variables.
The term ``constraint-coupled'' that we associate to problem~\eqref{eq:MILP}
is due to the fact that the constraints $\sum_{i=1}^N A_i \bx_i = b$ create
a link among all the variables $x_1, \ldots, x_N$, which otherwise could be
optimized independently from each other.
To achieve the mentioned reformulation, we now specify the quantities $\bx_i$,
$c_i$, $\Xi$, $A_i$ for each type of device, and the right-hand side vector $b$.

\emph{Storages}.
We assume that each storage unit $i \in \STOR$ is responsible for the
optimization vector $x_i$ consisting of the stack of $x_i(k), u_i(k), z_i(k) \in \real$
and $\delta_i(k) \in \{0, 1\}$ for all $k \in \fromto{0}{K-1}$ plus the variable
$x_i(K) \in \real$. The constraints in $\Xi$ are given by~\eqref{eq:storage},
while the cost function is
$c_i^\top x_i = \sum_{k=0}^{K-1} \OM_i (2 z_i(k) - u_i(k))$.

\emph{Generators}.
Each generator $i \in \GEN$ is responsible for the optimization vector
$x_i$ consisting of the stack of
$u_i(k), \nu_i(k), \theta_i^\textsc{u}(k), \theta_i^\textsc{d}(k) \in \real$
and $\delta_i(k) \in \{0, 1\}$ for all $k \in \fromto{0}{K-1}$.
The constraints in $\Xi$ are given by~\eqref{eq:generator_1}--\eqref{eq:generator_9},
while the cost function is
$c_i^\top \bx_i = \sum_{k=0}^{K-1} \big( \OM_i \delta_i(k) + \nu_i(k) +
\theta_i^\textsc{u}(k) + \theta_i^\textsc{d}(k) \big)$.

\emph{Critical loads}.
For the critical loads $i \in \LOAD$ there are no variables to optimize,
but they must be taken into account in the coupling constraints.

\emph{Controllable loads}.
For each controllable load $i \in \CLOAD$ the optimization vector
$\bx_i$ consists of the stack of $\beta_i(k) \in \real$,
for all $k \in \fromto{0}{K-1}$, with constraints given by~\eqref{eq:load}.
Note that, for this class of devices, the local
constraint set is not mixed-integer. The cost function is
$c_i^\top \bx_i = \sum_{k=0}^{K-1} \varphi_i D_i(k) \beta_i(k)$.

\emph{Connection to the utility grid}.
For this device $i \in \UGRID$, the optimization vector $\bx_i$ consists
of the stack of $u_i(k), \phi_i(k) \in \real$ and $\delta_i(k) \in \{0, 1\}$
for all $k \in \fromto{0}{K-1}$.
The local constraints are~\eqref{eq:grid}, while the cost function is
$c_i^\top \bx_i = \sum_{k=0}^{K-1} \phi_i(k)$.

\emph{Coupling constraints}.
Finally, the coupling constraints are given by~\eqref{eq:power_balance},
which can be encoded in the form $\sum_{i=1}^N A_i x_i = b$ by
appropriately defining the matrices $A_i$ and the vector $b$.
\begin{subequations}
In particular, the matrices $A_i \in \real^{K \times n_i}$ are such that
\begin{align}
  [A_i x_i]_k &= u_i(k), &&\forall i \in \STOR,
  \\
  [A_i x_i]_k &= -u_i(k), &&\forall i \in \GEN,
  \\
  [A_i x_i]_k &= -\beta_i(k) D_i(k), &&\forall i \in \CLOAD,
  \\
  [A_i x_i]_k &= -u_i(k), &&\forall i \in \UGRID,
\end{align}
for all times $k$, while the right-hand side vector $b \in \real^K$ is equal to
\begin{align}
  b = - \sum_{i \in \CLOAD} D_i - \sum_{i \in \LOAD} D_i + \sum_{i \in \REN} P_i,
\label{eq:h_def}
\end{align}%
\end{subequations}
where here $D_i \in \real^K$ and $P_i \in \real^K$ denote the
stack of $D_i(k)$ and $P_i(k)$ for all times $k$.
Note that the power generated by the renewables introduces a stochasticity
in the right-hand side vector $b$ appearing in problem~\eqref{eq:MILP}.

In the considered distributed context,
we assume that each agent $i$ does not know the entire problem information.
In particular, we assume it only knows the local cost vector $c_i$, the local
constraint $\Xi$ and its matrix $A_i$ of the coupling constraint.
The exchange of information among $N$ agents occurs according to a
graph-based
communication model. We use $\GG = (V, \EE)$ to indicate the undirected,
connected graph describing the network, where $V = \until{N}$ is the set of vertices
and $\EE$ is the set of edges.
If $(i,j)\in \EE$, then agent $i$ can communicate with agent $j$ and viceversa.
We use $\nbrs_i$ to indicate the set of neighbors of agent $i$ in $\GG$, i.e.,
$\nbrs_i = \{ j \in V | (i,j) \in \EE \}$.

\subsection{Two-stage Stochastic Optimization Approach}
In its current form, problem~\eqref{eq:MILP} cannot be practically solved due
to the right-hand side vector $b$ being unknown. To deal with this,
the approach consists
of considering a set of possible \emph{scenarios} that may arise
and then to formulate and solve a so-called \emph{two-stage}
stochastic optimization problem, which we now introduce.

Intuitively, in this uncertain scenario one has to ``a priori'' (i.e. without knowing the actual value of the random
vector $b$) choose a set of
control actions $u_i(k)$, such as generated/stored power or power curtailments,
in order to minimize a certain cost criterion in an expected sense. However, these
control actions will inevitably result in a violation of the power
balance constraint~\eqref{eq:power_balance}
``a posteriori'' (i.e. when the actual power production
of renewables, and hence value of the random vector $b$,
becomes known). To compensate for this infeasibility,
\emph{recourse} actions must be taken. These actions are associated to a
cost and will have an impact on the final performance achieved
by the whole control scheme.
In the jargon of two-stage stochastic optimization, the first-stage
optimization variables are those associated to the control actions
(i.e. $x_1, \ldots, x_N$ in problem~\eqref{eq:MILP}), while
the second-stage optimization variables (to be introduced shortly)
are those associated to recourses.

Formally, we denote by $\omega$ the random
vector collecting all the renewable energy generation profiles.
We assume a finite discrete probability distribution for
$\omega$ and we denote by $\pi_r$ the probability of each $\omega_r$,
i.e. $\pi_r = \mathbb{P} (\omega = \omega_i)$ for all $r \in \until{R}$.
To keep the notation consistent we denote the renewable energy profile
corresponding to $\omega_r$ as $P_{ir}(k)$.
We denote by $b_r$ the realization of $b$ associated to the scenario $\omega_r$.
Using these positions, the two-stage stochastic MILP can be
formulated as
\begin{align}
\begin{split}
  \min_{\substack{\bx_1,\ldots,\bx_N \\ \eta_+, \eta_-}} \:
  & \: \sum_{i =1}^N c_i^\top \bx_i + \sum_{k=0}^{K-1} \sum_{r=1}^R \pi_r (q_+ \eta_{kr+} + q_- \eta_{kr-})
  \\
  \subj \: 
  & \: -\eta_{r-} \le \sum_{i=1}^N A_i x_i - b_r \le \eta_{r+}, \hspace{0.5cm} r = \interv{R}
  \\
  & \: \eta_{+}, \eta_{-} \ge \0,
  \\
  & \: \bx_i \in \Xi, \hspace{1.5cm} i = \interv{N},
\end{split}
\label{eq:two_stage_MILP}
\end{align}
where $x_1, \ldots, x_N$ are the first-stage variables modeling the (a-priori) control actions
and $\eta_+, \eta_-$ are the two-stage variables modeling the (a-posteriori) recourse actions,
which are penalized in the cost with $q_+ \ge 0$ and $q_- \ge 0$, which are the costs related
to energy surplus and shortage, respectively.
In problem~\eqref{eq:two_stage_MILP}, we denoted by $\eta_{kr+}$ the variable associated with positive
recourse for scenario $r$ at time $k$.
We also use the symbol $\eta_{r+}$ to denote
the stack of $\eta_{kr+}$ for all $k$.
The stack of $\eta_{kr+}$ for all $k$ and $r$ is denoted
by $\eta_{+}$. A similar notation holds for $\eta_-$.
It can be seen that the additional term
in the cost is the expected value of the cost associated to recourse actions, i.e.
\begin{align*}
  &\sum_{k=0}^{K-1} \!\sum_{r=1}^R \!\pi_r (q_+ \eta_{kr+} \!+\! q_- \eta_{kr-})
  \!=\!\! \sum_{k=0}^{K-1} \!\expv \!\bigg[ \!\Phi \bigg( \!\bigg[\!\sum_{i=1}^N A_{i} x_i - b \bigg]_k \bigg) \!\bigg]\!,
\end{align*}
where $\Phi(z) = q_+ z$ if $z \ge 0$ and $\Phi(z) = - q_- z$ if $z < 0$.

At a first glance, it may seem that the two-stage problem~\eqref{eq:two_stage_MILP}
loses the constraint-coupled structure of the distributed optimization problem~\eqref{eq:MILP}.
However, with a bit a manipulation, it is still possible to arrive at a similar result.
We begin by streamlining the notation. Define $\eta \in \real^{2KR}$, $\eta \ge 0$ as the
stack of $\eta_{+}$ and $\eta_{-}$, and the vector $d \in \real^{2KR}$
such that $d^\top \eta = \sum_{k=0}^{K-1} \sum_{r=1}^R \pi_r (q_+ \eta_{kr+} + q_- \eta_{kr-})$.
Moreover, define $H_i \in \real^{2KR \times n_i}$ and $h \in \real^{2KR}$ with
\begin{align*}
  H_i &= \1 \kron \begin{bmatrix} A_i^\top \!&\! -A_i^\top \end{bmatrix}^\top
   \!= \begin{bmatrix}
    A_i^\top \!&\!
    -A_i^\top \!&\!
    \cdots
    \!&\!
    A_i^\top \!&\!
    -A_i^\top
  \end{bmatrix}^\top\!,
  \\
  h &= \begin{bmatrix}
    b_{1}^\top &
    -b_{1}^\top &
    \cdots
    &
    b_{R}^\top &
    -b_{R}^\top
  \end{bmatrix}^\top,
\end{align*}
where $\1 \in \real^R$ is the vector of ones and $\kron$ denotes the kronecker product.
Thus, problem~\eqref{eq:two_stage_MILP} is equivalent to
\begin{align}
\begin{split}
  \min_{\substack{\bx_1,\ldots,\bx_N \\ \eta}} \:
  & \: \sum_{i =1}^N c_i^\top \bx_i + d^\top \eta
  \\
  \subj \: 
  & \: \sum_{i=1}^N H_i x_i - h \le \eta
  \\
  & \: \eta \ge \0,  \:\:
  \bx_i \in \Xi, \hspace{1cm} i = \interv{N}.
\end{split}
\label{eq:two_stage_MILP_stream}
\end{align}
By defining $\eta_1, \ldots, \eta_N \in \real^{2RK}$
such that $\sum_{i=1}^N \eta_i = \eta$ and each $\eta_i \ge 0$,
we see that problem~\eqref{eq:two_stage_MILP_stream}
is finally equivalent to
\begin{align}
\begin{split}
  \min_{\substack{\bx_1,\ldots,\bx_N \\ \eta_1, \ldots, \eta_N}} \:
  & \: \sum_{i =1}^N (c_i^\top \bx_i + d^\top \eta_i)
  \\
  \subj \: 
  & \: \sum_{i=1}^N (H_i x_i - \eta_i) \le h,
  \\
  & \: \eta_i \ge \0, \:\: \bx_i \in \Xi, \hspace{1cm} i = \interv{N},
\end{split}
\label{eq:two_stage_MILP_distr}
\end{align}
in the sense that any solution of~\eqref{eq:two_stage_MILP_stream}
can be reconstructed from a solution of~\eqref{eq:two_stage_MILP_distr} by using
$\eta = \sum_{i=1}^N \eta_i$.
Note that problem~\eqref{eq:two_stage_MILP_distr} has an unbounded feasible
set (because of the variables $\eta_i$) but it always admits an optimal solution
due to the terms $d^\top \eta_i$ minimized in the cost
(recall that $d \ge \0$).

\section{Distributed Algorithm and Analysis}
\label{sec:algorithm}

We now propose a distributed algorithm to compute a
feasible solution to problem~\eqref{eq:two_stage_MILP_distr}
and provide the convergence results.

\subsection{Distributed Algorithm Description}
Let us begin by describing the proposed distributed algorithm to solve
problem~\eqref{eq:two_stage_MILP_distr}.
The basic idea behind the distributed algorithm is
to compute a mixed-integer solution starting from an optimal
solution of the convex relaxation of problem~\eqref{eq:two_stage_MILP_stream}
obtained by replacing $\Xi$ with their convex hull $\convXi$,
\begin{align}
\begin{split}
  \min_{\substack{\bz_1,\ldots,\bz_N \\ \eta_1, \ldots, \eta_N}} \:
  & \: \sum_{i =1}^N (c_i^\top \bz_i + d^\top \eta_i)
  \\
  \subj \: 
  & \: \sum_{i=1}^N (H_i \bz_i - \eta_i) \le h,
  \\
  & \: \eta_i \ge \0, \:\: \bz_i \in \conv{\Xi}, \hspace{0.5cm} i = \interv{N},
\end{split}
\label{eq:two_stage_MILP_stream_approx}
\end{align}
where we denote by $\bz_i$ the continuous counterpart of the mixed-integer variable $\bx_i$.
To do so, each agent $i$ maintains an auxiliary variable
$\by_i^t \in \real^{2RK}$, which represents a local \emph{allocation}
of the coupling constraints (cf. Appendix~\ref{sec:primal_decomp}).
At each iteration $t$, the vector $\by_i^t$
is updated according to~\eqref{eq:alg_z_LP}--\eqref{eq:alg_y_update}.
After $T_f > 0$ iterations, the agent computes
a tentative mixed-integer solution based
on the last computed allocation estimate (cf.~\eqref{eq:alg_x_MILP}).
Algorithm~\ref{alg:sg_algorithm} summarizes the steps from the
perspective of agent $i$.
\begin{algorithm}[htbp]
\floatname{algorithm}{Algorithm}

  \begin{algorithmic}[0]
  
    \Statex \textbf{Initialization}:
      set $T_f > 0$ and
      $\by_{i}^0$ such that $\sum_{i=1}^N \by_i^0 = h$
    \smallskip

    \Statex \textbf{Repeat} for $t = 0, 1, \ldots, T_f-1$ %
    \smallskip
    
      \StatexIndent[0.75]
      \textbf{Compute} $\bmu_i^t$ as a Lagrange multiplier of %
      \begin{align} 
      \label{eq:alg_z_LP}
      \begin{split}
        \min_{\bz_{i}, \eta_i } \hspace{1.2cm} &\: c_i^\top \bz_{i} + d^\top \eta_i
        \\
        \subj \hspace{0.3cm} 
        \: \bmu_i : \: & \: H_i \bz_i \leq \by_i^t + \eta_i
        \\
        & \: \eta_i \ge \0, \:\: \bz_i \in \conv{\Xi}
      \end{split}
      \end{align}

      \StatexIndent[0.75]
      \textbf{Receive} $\bmu_{j}^t$ from $j\in\nbrs_i$ and update %
      \begin{align}
      \label{eq:alg_y_update}
      \begin{split}
        \by_{i}^{t+1} = \by_{i}^t + \alpha^t \sum_{j \in \nbrs_i} \big( \bmu_{i}^t - \bmu_{j}^t \big)
      \end{split}
      \end{align}

      \Statex
      \textbf{Return} $(\bx_i^{T_f}, \eta_i^{T_f})$ as optimal solution of
      \begin{align}
      \label{eq:alg_x_MILP}
      \begin{split}
        \min_{\bx_{i}, \eta_i } \: &\: c_i^\top \bx_{i} + d^\top \eta_i
        \\
        \subj \: & \: H_i \bx_i \leq \by_i^{T_f} + \eta_i
        \\
        & \: \eta_i \ge \0, \:\: \bx_i \in \Xi
      \end{split}
      \end{align}

  \end{algorithmic}
  \caption{\sgalgname/}
  \label{alg:sg_algorithm}
\end{algorithm}

Let us briefly comment on the algorithm structure.
As it will be clear from the analysis,
the first two steps~\eqref{eq:alg_z_LP}--\eqref{eq:alg_y_update}
are used to compute an optimal solution of problem~\eqref{eq:two_stage_MILP_stream_approx},
while the last step~\eqref{eq:alg_x_MILP} reconstructs a mixed-integer solution.
Note that problem~\eqref{eq:alg_z_LP} is an LP and
problem~\eqref{eq:alg_x_MILP} is a MILP. From a computational
point of view, in order to compute a Lagrange multiplier of
problem~\eqref{eq:alg_z_LP} the agent can locally run either a dual
subgradient method or a dual cutting-plane method (cf.~\cite{camisa2021distributed}),
while an optimal solution to problem~\eqref{eq:alg_x_MILP} can be
found with any MILP solver.
In the next subsection we will prove
a worst-case violation of the power balance constraints.

\begin{remark}
	An important fact is that the computed mixed-integer solution
	always satisfies the coupling constraint appearing in
	problem~\eqref{eq:two_stage_MILP_stream}
  with a possibly high $\eta_i$, i.e.
	\begin{align*}
	  \sum_{i=1}^N (H_i x_i^{T_f} - \eta_i^{T_f})
	  \le
	  \sum_{i=1}^N y_i^{T_f}
	  =
	  h,
	\end{align*}
	where the inequality follows by construction and the equality
	follows by the forthcoming Lemma~\ref{lemma:DPD_convergence_sg}.
	Thus, the algorithm can be stopped at any iteration $T_f \ge 0$ and the
	resulting solution will be feasible for the two-stage
	MILP~\eqref{eq:two_stage_MILP_stream}.
	The greater the number of iterations, the higher is the optimality
	of the computed solution and the lower is the expected violation
	of the original power balance constraint.
  \oprocend
\end{remark}

\subsection{Theoretical Results}
In this subsection, we provide theoretical results on Algorithm~\ref{alg:sg_algorithm}.
In particular, we will prove a bound for the
worst-case violation of the asymptotically computed mixed-integer solution.
To begin with, we recall some preliminary lemmas, where we remind that $K$
denotes the prediction horizon and $R$ is the total number of scenarios in
the stochastic problem).
\begin{lemma}[\cite{camisa2021distributed}]
\label{lemma:LP_integer_components_sg}
	Let problem~\eqref{eq:two_stage_MILP_stream_approx} be feasible and let
	$(\bar{\bz}_1, \ldots, \bar{\bz}_N, \bar{\eta}_1, \ldots, \bar{\eta}_N)$
	be any vertex of its feasible set.
	Then, there exists an index set $\INT \subseteq \until{N}$,
	with cardinality $|\INT| \ge N-2RK$,
	such that $\bar{\bz}_i \in \Xi$ for all $i \in \INT$.
	\oprocend
\end{lemma}
The consequence of Lemma~\ref{lemma:LP_integer_components_sg} is that
at least $N-2RK$ blocks of the mixed-integer solution computed
asymptotically by Algorithm~\ref{alg:sg_algorithm} are equal to the corresponding
blocks of optimal solution of~\eqref{eq:two_stage_MILP_stream_approx}.
Next we recall convergence of the steps~\eqref{eq:alg_z_LP}--\eqref{eq:alg_y_update}.
To this end, we denote as $(\bz_1^\LP, \ldots, \bz_N^\LP, \eta_1^\LP, \ldots, \eta_N^\LP)$
an optimal solution of problem~\eqref{eq:two_stage_MILP_stream_approx}, together with
the allocation vector $(y_1^\LP, \ldots, y_N^\LP)$ associated to the primal decomposition
master problem (cf. Appendix~\ref{sec:primal_decomp}),
which is a vector satisfying
\begin{subequations}
\begin{align}
  H_i z_i^\LP - \eta_i^\LP &\le y_i^\LP, \hspace{0.3cm} \text{ for all } i \in \until{N},
  \\
  \text{and} \hspace{0.3cm}
  \sum_{i=1}^N y_i^\LP &= h.
\end{align}\end{subequations}
The following assumption is made on the step-size sequence.
\begin{assumption}
  \label{ass:step-size_sg}
  The step-size sequence $\{ \alpha^t \}_{t\ge0}$, with each $\alpha^t \ge 0$,
  satisfies $\sum_{t=0}^{\infty} \alpha^t = \infty$,
  $\sum_{t=0}^{\infty} \big( \alpha^t \big)^2 < \infty$.
  \oprocend
\end{assumption}
The following proposition summarizes the convergence
properties of the steps~\eqref{eq:alg_z_LP}--\eqref{eq:alg_y_update}.
\begin{lemma}[\cite{camisa2021distributed}]
\label{lemma:DPD_convergence_sg}
	Let problem~\eqref{eq:two_stage_MILP_stream_approx} be feasible and let
	Assumption~\ref{ass:step-size_sg} hold. Consider the allocation vector sequence
	$\{ \by_1^t,\ldots,\by_N^t \}_{t\ge 0}$ generated by
	steps~\eqref{eq:alg_z_LP}--\eqref{eq:alg_y_update}
	of Algorithm~\ref{alg:sg_algorithm}
	with the allocation vectors $\by_i^0$ initialized such that
	$\sum_{i=1}^N \by_i^0 = h$.
	Then,
  \begin{enumerate}
	  \item $\sum_{i=1}^N \by_i^t  = h$ for all $t \ge 0$;
	   
	  \item $\lim_{t \to \infty} \| \by_i^t - \by_i^\LP \| = 0$ for all
	    $i \in \until{N}$.
	  \oprocend
  \end{enumerate}
\end{lemma}
Because of Lemma~\ref{lemma:DPD_convergence_sg}, from now
on we concentrate on the asymptotic mixed-integer solution computed
by Algorithm~\ref{alg:sg_algorithm}. In particular, we denote by
$(\algx_i, \algeta_i)$ the optimal solution of problem~\eqref{eq:alg_x_MILP}
with allocation equal to $y_i^\LP$, i.e.
\begin{align}
\label{eq:alg_x_MILP_asymptotic}
\begin{split}
  \min_{\bx_{i}, \eta_i } \: &\: c_i^\top \bx_{i} + d^\top \eta_i
  \\
  \subj \: & \: H_i \bx_i \leq \by_i^\LP + \eta_i
  \\
  & \: \eta_i \ge \0, \:\: \bx_i \in \Xi.
\end{split}
\end{align}
We also define the lower bound of resources $\bL_i \in \real^{2RK}$
\begin{align*}
  \bL_i \triangleq \min_{x_i, \eta_i} \: & \: H_i x_i - \eta_i
  \\
  \subj \: & \: x_i \in \conv{X_i}
  \\
  & \: \0 \le \eta_i \le M \1.
\end{align*}
where $\min$ is component-wise and $M>0$ is a sufficiently large number.
Thus, it holds $\ell_i \le y_i$ for all
admissible allocations $y_i$, and in particular $\ell_i \le y_i^\LP$.
Operatively, since the constraints on $x_i$
and $\eta_i$ are disjoint the vector $\bL_i$ can be computed by replacing
$x_i \in \conv{X_i}$ with $x_i \in X_i$.
In the next theorem we formalize the bound on the worst-case violation.
\begin{theorem}
\label{thm:worst_case_viol}
  Let problem~\eqref{eq:two_stage_MILP_stream_approx} be feasible
  and consider the asymptotic mixed-integer solution $(\algx_i, \algeta_i)$
  computed by each agent $i \in \until{N}$.
  Then, the worst-case violation of the power balance constraint is
  \begin{align*}
	  \sum_{i=1}^N H_i \algx_i - h
	  \le
	  \sum_{i \in \INT} \eta_i^\LP + \sum_{i \notin \INT} \frac{c_i^\top (\xl_i - \algx_i) + d^\top \etal_i}{d^\textsc{min}} \1,
	\end{align*}
	where $d^\textsc{min} = \min_{j \in \until{2RK}} d_j$, $\INT$ denotes
	the set of agents (satisfying $|\INT| \ge N-2RK$) for which $z_i^\LP \in \Xi$
	and $(\xl_i, \etal_i)$ is an optimal solution of problem~\eqref{eq:proof_worst_case_xil}.
	\oprocend
\end{theorem}
The proof is provided in Appendix~\ref{sec:proof_theorem}.
Note that, since this bound is the sum of contributions of
the agents, it can be computed a posteriori in a distributed way
using a consensus scheme. To do so, they first need to detect
whether they belong to $\INT$ or not by computing the primal
solution $z_i^\LP$ of~\eqref{eq:alg_z_LP} and by checking
whether it satisfies $z_i^\LP \in \Xi$. Then, they run the
consensus scheme using as initial condition either $N \eta_i^\LP$
(if $z_i^\LP \in \Xi$) or $N\frac{c_i^\top (\xl_i - \algx_i) + d^\top \etal_i}{d^\textsc{min}} \1$
(if $z_i^\LP \notin \Xi$).

\section{Numerical Experiments}
\label{sec:simulations}

In this section, we validate the proposed framework through large-scale
numerical computations. All the simulations are performed with the \textsc{disropt}
package~\cite{farina2019disropt} and are performed on the Italian
HPC CINECA infrastructure. In order to make the simulations
realistic, we run Algorithm~\ref{alg:sg_algorithm} on a generated problem
with data synthesized using a deep Generative Adversarial Network
(GAN)~\cite{goodfellow2014generative}. In the next subsections, we
first provide details regarding the scenario generation for renewable
energy sources, then we show aggregate results on Monte Carlo simulations
and finally we show in more detail one specific simulation.

\subsection{Scenario Generation with Generative Adversarial Networks}
\label{sec:gan}
Recall that $b \in \real^{K}$ is a random variable that
depends on the total energy produced by the renewables~\eqref{eq:h_def}.
The variable $b$ has its own probability distribution and
$b_1, \ldots, b_R \in \real^K$ are
randomly drawn samples (cf.~\eqref{eq:two_stage_MILP}).
In order to generate such samples, we utilize a Generative Adversarial Network
trained with an open historical dataset from the EU. To train the neural network, we
used the data series provided by Open Power System Data \cite{opsd2020timeseries}.
In particular, we used the generation data of renewable energy sources
in South Italy. To guarantee a certain uniformity of the data, we narrowed the
dataset by concentrating only on summer months
and discarded days with missing information.
Each sample is a vector in
$\real^{24}$ and contains information on
the power produced
during a day with a hourly resolution.

As for the utilized neural networks, the generative networks have
a $10$-dimensional input with the following layers:
\begin{itemize}
  \item a dense layer with $1536$ units, batch normalization and Leaky ReLU activation function;
  
  \item a layer that reshapes the input to the shape $(6, 256)$;
  
  \item a transposed convolution layer with $128$ output filters, kernel size equal to $5$, stride $1$,
    batch normalization and Leaky ReLU activation function;
  
  \item a transposed convolution layer with $64$ output filters, kernel size equal to $5$, stride $2$,
    batch normalization and Leaky ReLU activation function;
  
  \item a transposed convolution layer with $1$ output filter, kernel size equal to $5$, stride $2$
    and $\tanh$ activation function.
\end{itemize}
The ouput of the generative network is a $24$-dimensional vector containing
the power produced by the renewable unit at each time slot of the day.
The discriminator networks have a $24$-dimensional input with the following
layers:
\begin{itemize}
  \item a convolution layer with $64$ output filters, kernel size equal to $5$, stride $2$ and
    Leaky ReLU activation function;
  
  \item a Dropout layer with rate $0.3$;
  
  \item a convolution layer with $128$ output filters, kernel size equal to $5$, stride $2$ and
    Leaky ReLU activation function;
  
  \item a Dropout layer with rate $0.3$;
  
  \item a layer that flattens the input;
  
  \item a dense layer with one output unit.
\end{itemize}
The output of the discriminator networks is a scalar that denotes the
probability that the evaluated input is a real one or a generated one.

\begin{figure}[htbp]\centering
  \includegraphics[scale=1]{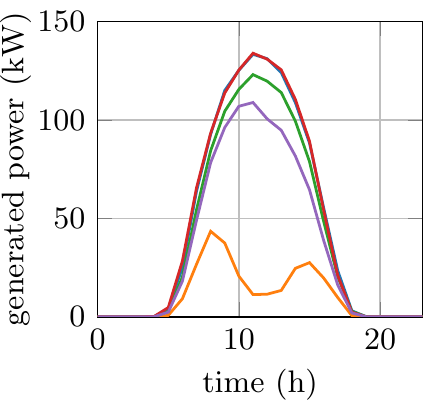}
  \hfill
  \includegraphics[scale=1]{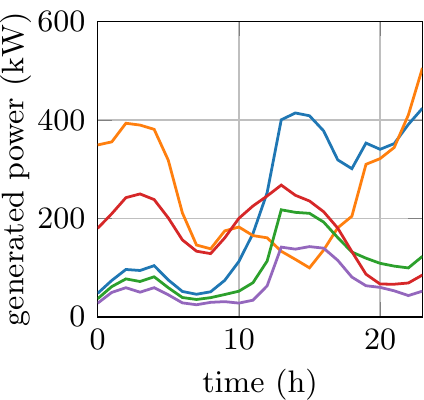}
  
  \caption{
    Five examples of power generation profiles generated by the GANs.
    Left: solar energy, right: wind energy.
  }
\label{fig:gan}
\end{figure}

We used the neural networks to generate samples of solar energy and
wind energy. We used \textsc{Tensorflow} 2.4 to model the networks
and we performed the training with $10^4$ epochs using the ADAM
algorithm.
In Figure~\ref{fig:gan}, we show example profiles of
solar and wind energy generated by the networks.
It can be noted that generated trajectories of solar energy production have a
maximum at midday, while one of the trajectories has lower values than the
others and may be associated, for instance, to a cloudy day. In any case,
the power generated outside the time window 5am-8pm is close to zero,
consistently with real profiles.

\subsection{Monte Carlo Simulations}
To test the proposed framework, we performed $100$
Monte Carlo simulations in which we run
Algorithm~\ref{alg:sg_algorithm} on different realizations
of the energy generation scenarios (i.e., different realizations of $b$).

\begin{figure}[htbp]\centering
  \includegraphics[scale=1]{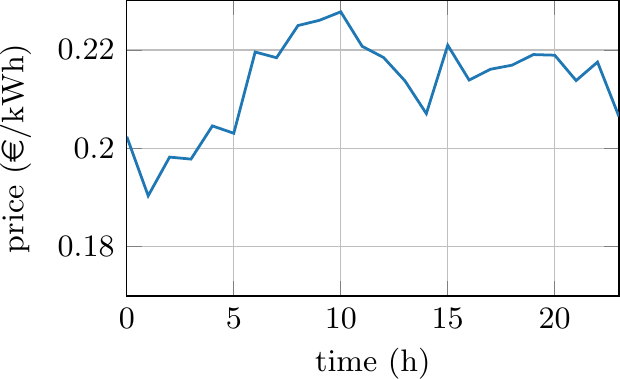}
  \caption{
    Daily spot prices from Open Power System Data \cite{opsd2020timeseries}.
  }
\label{fig:stochastic_daily_spot_prices}
\end{figure}

We considered a microgrid control problem with the following units: $20$ generators,
$20$ storages, $60$ controllable loads, $20$ critical loads, $40$ solar generators,
$15$ wind generators and the connection to the main grid.
For each instance of the problem, we extracted $R = 5$ scenarios and fixed a
$24$-hour prediction horizon and $1$-hour sampling time.
The initial conditions of storages and generators are generated randomly.
As regards the load profiles and the daily spot prices,
we utilized the data provided by~\cite{opsd2020timeseries}, which
are shown in Figure~\ref{fig:stochastic_daily_spot_prices}.
We then executed Algorithm~\ref{alg:sg_algorithm} for $500$ iterations with
a piece-wise constant step size that we initialize to $3.0$ and multiply by
$0.5$ every $100$ iterations.

\begin{figure}[htbp]\centering
  \includegraphics[scale=1]{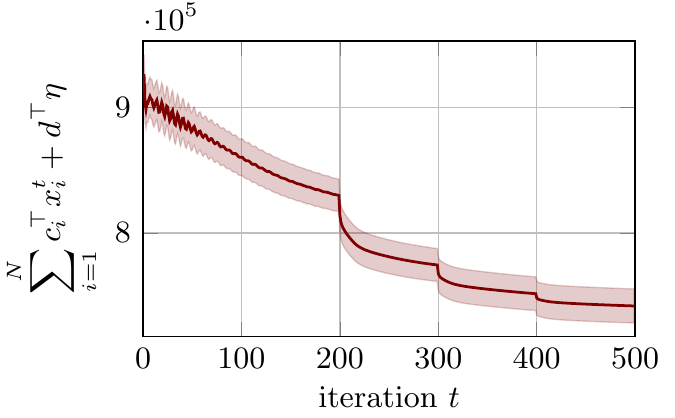}
  \caption{
    Evolution of the cost yielded by Algorithm~\ref{alg:sg_algorithm}. The solid
    line is the mean value of the Monte Carlo trials while the dashed area represents
    one standard deviation.
  }
\label{fig:montecarlo_cost}
\end{figure}

\begin{figure}[htbp]\centering
  \includegraphics[scale=1]{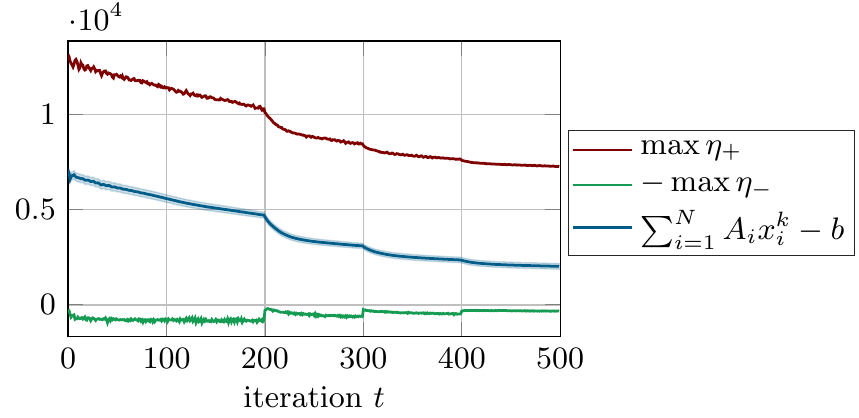}
  \caption{
    Evolution of the coupling constraint value throughout the
    evolution of Algorithm~\ref{alg:sg_algorithm}. The blue line represents
    the average value of the power balance constraint (the dashed area
    corresponds to one standard deviation). The upper and lower lines
    are the maximum positive and negative two-stage violations
    of the constraints.
  }
\label{fig:montecarlo_coupling}
\end{figure}

The results of the simulations are shown in Figures~\ref{fig:montecarlo_cost}
and~\ref{fig:montecarlo_coupling}. In Figure~\ref{fig:montecarlo_cost},
we plot the cost of the mixed-integer solution computed by the algorithm
throughout its evolution (in particular the cost function of the two-stage
problem~\eqref{eq:two_stage_MILP}). The picture highlights how
the algorithm improves the cost at each iteration, i.e., the more iterations
are performed, the more the solution performance improves.

In Figure~\ref{fig:montecarlo_coupling}, we show the value of the coupling
constraints for the two-stage problem~\eqref{eq:two_stage_MILP}.
The red and green lines correspond to the maximum value of $\eta_+$
and $\eta_-$ (with changed sign) respectively, where the maximum
is taken with respect to the scenarios, to the components
of the constraint and to the Monte Carlo trials.
The blue line represents the average value of the
power balance constraint, while the dashed area corresponds to one
standard deviation of the Monte Carlo trials.
At each time step the power balance constraints are always
in between the upper and lower line, while the uncertainty range
reduces as the algorithm progresses.

\subsection{Results on a Single Instance}
\label{sec:simulation_stochastic}
\begin{figure}[htbp]\centering
  \includegraphics[scale=1]{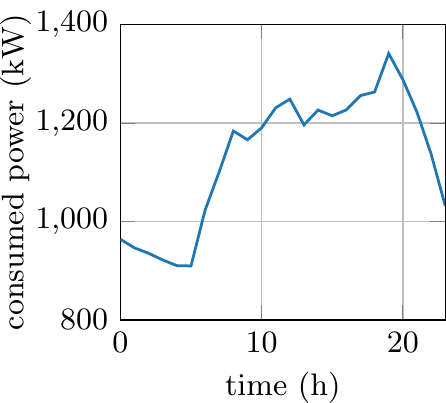}
  \hfill
  \includegraphics[scale=1]{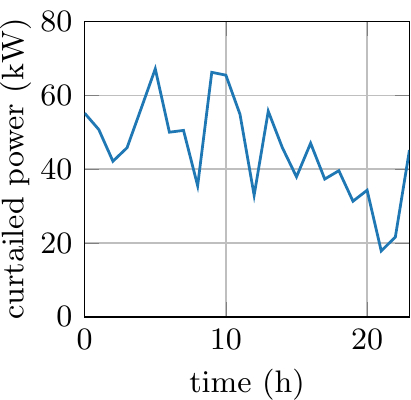}
  
  \caption{
    Total consumed power (critical and controllable loads) and curtailed power (for controllable loads only).
  }
\label{fig:stochastic_consumed_curtailed_power}
\end{figure}

\begin{figure}[htbp]\centering
  \includegraphics[scale=1]{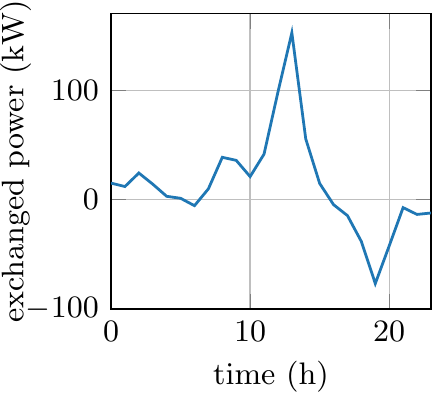}
  \hfill
  \includegraphics[scale=1]{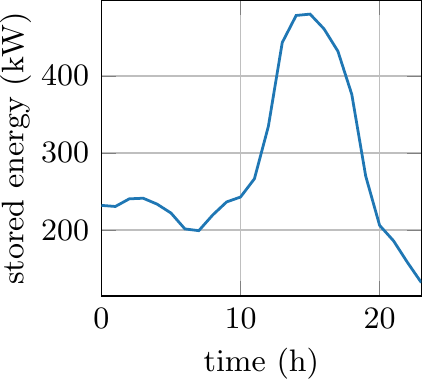}
  
  \caption{
    Total average power exchanged by storage units (left) and level of total stored power (right).
  }
\label{fig:stochastic_storages}
\end{figure}

\begin{figure}[htbp]\centering
  \includegraphics[scale=1]{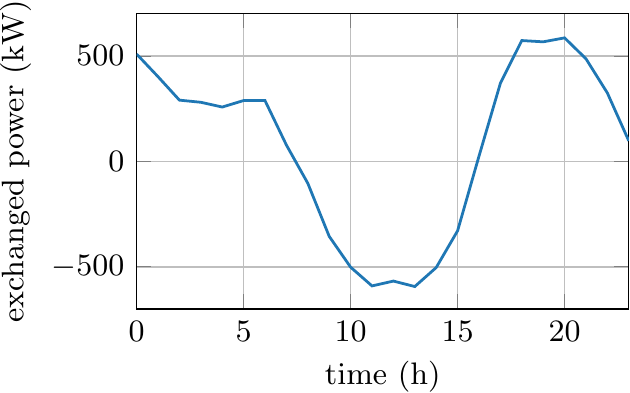}
  \caption{
    Total power exchanged with the utility grid.
  }
\label{fig:stochastic_grid}
\end{figure}

\begin{figure}[htbp]\centering
  \includegraphics[scale=1]{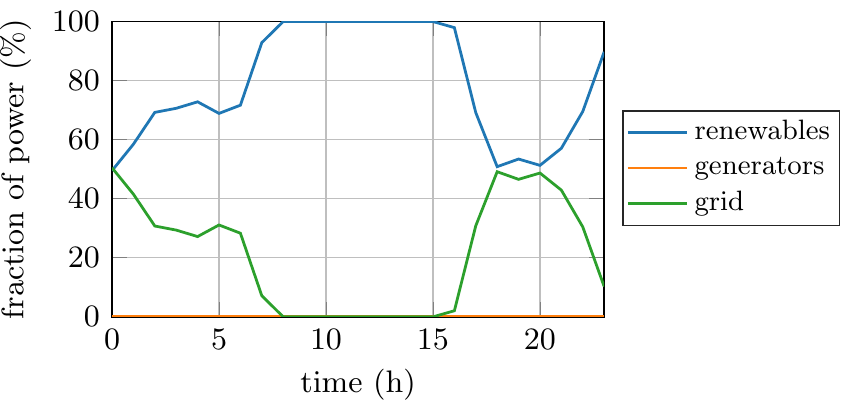}
  \caption{
    Fraction of consumed power coming from generators, renewables and utility grid at each time slot.
  }
\label{fig:stochastic_power_fraction}
\end{figure}

To conclude this section, we show how Algorithm~\ref{alg:sg_algorithm}
behaves on a single instance of the Monte Carlo trials.
In Figure~\ref{fig:stochastic_consumed_curtailed_power}, we show the
total consumed power and the total curtailed power.
In Figure~\ref{fig:stochastic_storages}, we show the total power exchanged
with storage units (a positive value means that, overall, the storage units are
charging) and the global level of stored power.
The solution provided by the algorithm is such that storages
accumulate as much energy as they can during the peaks of power produced
by the renewables. This energy is then released during the subsequent hours
of the day.
In Figure~\ref{fig:stochastic_grid}, we show the total power exchanged
with the utility grid (a positive value means that power is purchased
from the grid).
Note that, during the peak of power produced by the renewables,
the microgrid exports energy to the main grid
in order to maximize the income.
In Figure~\ref{fig:stochastic_power_fraction}, we show where does the
total available power comes from. In particular we highlight the
fraction of power coming from generators, renewables and the utility grid.
In this simulation, the generators did not produce any energy.

\section{Conclusions}
In this paper, we considered a microgrid control problem
to be solved over a peer-to-peer network of agents.
Each agent represents a unit of the microgrid and must cooperate
with the other units in order to solve the problem
without a centralized coordinator. We used a challenging stochastic
mixed-integer microgrid model and proposed a distributed
algorithm to solve the problem, for which we provided theoretical
guarantees on the constraint violation. Numerical computations
on a synthesized problem using Generative Adversarial Networks
show the validity of the proposed approach.

\appendix

\subsection{Review of Primal Decomposition}
\label{sec:primal_decomp}
Consider a network of $N$ agents indexed by $\agents = \until{N}$
that aim to solve a linear program of the form
\begin{align}
\begin{split}
  \min_{x_1, \ldots, x_N} \: & \: \sum_{i=1}^N c_i^\top x_i
  \\
  \subj \: & \: x_i \in X_i, \hspace{1cm} \forall i \in \agents,
  \\
  & \: \sum_{i=1}^N A_i x_i \le b,
\end{split}
\label{eq:app_LP}
\end{align}
where each $x_i \in \real^{n_i}$ is the $i$-th optimization variable,
$c_i \in \real^{n_i}$ is the $i$-th cost vector, $X_i \subset \real^{n_i}$
is the $i$-th polyhedral constraint set and $A_i \in \real^{S \times n_i}$
is a matrix for the $i$-th contribution to the \emph{coupling constraint}
$\sum_{i=1}^N A_i x_i \le b \in \real^S$. Problem~\eqref{eq:app_LP}
enjoys the constraint-coupled structure~\cite{notarstefano2019distributed}
and can be recast into a master-subproblem architecture
by using the so-called \emph{primal decomposition} technique~\cite{silverman1972primal}.
The right-hand side vector $b$ of the coupling constraint
is interpreted as a given (limited) resource to be shared among the
network agents.
Thus, local \emph{allocation vectors} $y_i \in \real^S$ for all $i$
are introduced such that $\sum_{i=1}^N y_i = b$.
To determine the allocations, a \emph{master problem} is introduced
\begin{align}
\begin{split}
  \min_{y_1,\ldots,y_N} \: & \: \sum_{i =1}^N p_i (y_i) 
  \\
  \subj \: & \: \smallsum_{i=1}^N y_i = b
  \\
  & \: y_i \in Y_i, \hspace{1cm} \forall i \in\agents,
\end{split}
\label{eq:app_primal_decomp_master}
\end{align}
where, for each $i\in\agents$, the function $\map{p_i}{\real^S}{\real}$ is
defined as the optimal cost of the $i$-th (linear programming) \emph{subproblem}
\begin{align}
\begin{split}
  p_i(y_i) = \: \min_{x_i} \: & \: c_i^\top x_i
  \\
  \subj \: 
  & \: A_i x_i \leq y_i
  \\
  & \: x_i \in X_i.
\end{split}
\label{eq:app_primal_decomp_subproblem}
\end{align}
In problem~\eqref{eq:app_primal_decomp_master}, the new constraint set
$Y_i \subseteq\real^S$ is
the set of $y_i$ for which
problem~\eqref{eq:app_primal_decomp_subproblem} is feasible, i.e.,
such that there exists $x_i \in X_i$ satisfying the local
\emph{allocation constraint} $A_i x_i \le y_i$.
Assuming problem~\eqref{eq:app_LP} is feasible and $X_i$ are compact sets,
if $(y_1^\star, \ldots, y_N^\star)$ is an optimal solution
of~\eqref{eq:app_primal_decomp_master} and, for all $i$,
$x_i^\star$ is optimal for~\eqref{eq:app_primal_decomp_subproblem}
(with $y_i = y_i^\star$), then $(x_1^\star, \ldots, x_N^\star)$
is an optimal solution of the original problem~\eqref{eq:app_LP}
(see, e.g.,~\cite[Lemma 1]{silverman1972primal}).

\subsection{Proof of Theorem~\ref{thm:worst_case_viol}}
\label{sec:proof_theorem}
	By the optimality of $(\algx_i, \algeta_i)$ for problem~\eqref{eq:alg_x_MILP_asymptotic}, it holds
	\begin{align}
	  c_i^\top \algx_i + d^\top \algeta_i \le c_i^\top x_i + d^\top \eta_i
	\label{eq:proof_worst_case_opt}
	\end{align}
	for all $x_i \in X_i$ and $\eta_i \ge 0$ such that $A_i \algx_i \le y_i^\LP + \eta_i$.
	One vector satisfying such condition is $(\xl_i, \etal_i)$ optimal solution of
	\begin{align}
	\begin{split}
	  \min_{x_i, \eta_i} \: & \: c_i^\top x_i + d^\top \eta_i
	  \\
	  \subj \: & \: \0 \le \eta_i \le M \1, \:\: x_i \in X_i,
	  \\
	  & \: H_i x_i \le \ell_i + \eta_i,
	\end{split}
	\label{eq:proof_worst_case_xil}
	\end{align}
	Indeed, it holds $H_i \xl_i \le \ell_i + \etal_i \le y_i^\LP + \etal_i$,
	where the first inequality is by construction and the second one
	follows by the discussion above on $\ell_i$.
	Thus, by using~\eqref{eq:proof_worst_case_opt} we conclude that
	\begin{align}
	  d^\top \algeta_i \le c_i^\top (\xl_i - \algx_i) + d^\top \etal_i.
	\end{align}
	By explicitly writing the scalar product $d^\top \algeta_i$ and by using
	the fact that $d, \eta \ge \0$ we obtain
	\begin{align*}
	  d^\top \algeta_i
	  = \sum_{j=1}^{2RK} d_j \algeta_{ij}
	  \ge \underbrace{\Big( \min_{j \in \until{2RK}} d_j \Big)}_{d^\textsc{min}} \sum_{j=1}^{2RK} \algeta_{ij}.
	\end{align*}
	Moreover, by using the fact that $\algeta_{ij} \le \sum_{k=1}^{2RK} \algeta_{ik}$
	for all $k$ we obtain
	\begin{align*}
	  \algeta_{i}
	  \le \frac{d^\top \algeta_i}{d^\textsc{min}} \1
	  \le \frac{c_i^\top (\xl_i - \algx_i) + d^\top \etal_i}{d^\textsc{min}} \1.
	\end{align*}
	Let us now compute an upper bound of the coupling constraint value, i.e.
	\begin{align}
	  \sum_{i=1}^N H_i \algx_i - h
	  &\le \underbrace{\sum_{i=1}^N y_i^\LP}_{b} + \sum_{i \in \INT} \eta_i^\LP + \sum_{i \notin \INT} \algeta_i - h
	  \nonumber
	  \\
	  &= \sum_{i \in \INT} \eta_i^\LP + \sum_{i \notin \INT} \algeta_i.
	\end{align}
	where we used the fact that, by Lemma~\ref{lemma:LP_integer_components_sg},
	for $i \in \INT$ it holds $H_i \algx_i \le y_i^\LP + \eta_i^\LP$,
	while for $i \notin \INT$ it holds $H_i \algx_i \le y_i^\LP + \algeta_i$.
	Thus we finally obtain the bound
	\begin{align*}
	  \sum_{i=1}^N H_i \algx_i - h
	  &\le
	  \sum_{i \in \INT} \eta_i^\LP + \sum_{i \notin \INT} \frac{c_i^\top (\xl_i - \algx_i) + d^\top \etal_i}{d^\textsc{min}} \1.
	\end{align*}
	and the proof follows.
	\oprocend

\bibliographystyle{IEEEtran}
\bibliography{smart_grid_biblio}

\end{document}